\documentstyle[righttag]{amsart}

\numberwithin{equation}{section}

\newcommand{\ms}{\medskip}

\newcommand{\noi}{\noindent}
\newcommand{\ra}{\rightarrow}
\newcommand{\bea}{\begin{eqnarray}}
\newcommand{\eea}{\end{eqnarray}}

\newcommand{\gr}{Groenewold}
\newcommand{\vh}{Van Hove}
\newcommand{\vn}{von Neumann}

\newcommand{\shs}{spherical harmonics}
\newcommand{\cg}{Clebsch-Gordan coefficient}
\newcommand{\ci}{S_i}
\newcommand{\cii}{S_i\,\!^2}
\newcommand{\ciii}{S_i\,\!^3}
\newcommand{\cj}{S_\ell}

\newcommand{\cx}{S_1}
\newcommand{\cxx}{S_1\,\!^2}

\newcommand{\cy}{S_2}
\newcommand{\cyy}{S_2\,\!^2}

\newcommand{\cz}{S_3}
\newcommand{\czz}{S_3\,\!^2}
\newcommand{\czzz}{S_3\,\!^3}
\newcommand{\q}{\cal Q}
\newcommand{\qx}{\cal Q(S_1)}
\newcommand{\qy}{\cal Q(S_2)}
\newcommand{\qz}{\cal Q(S_3)}
\newcommand{\qi}{\cal Q(S_i)}
\newcommand{\qii}{\cal Q\big(S_i\,\!^2\big)}
\newcommand{\qiii}{\cal Q\big(S_i\,\!^3\big)}
\newcommand{\qj}{\cal Q(S_\ell)}

\newcommand{\p}{\cal P}
\newcommand{\h}{\cal H}
\newcommand{\oo}{\cal O}
\newcommand{\ld}{\lambda}

\newcommand{\sm}{\Bbb S_{m}^{|m|}}

\newtheorem{thm}{Theorem}
\newtheorem{lem}{Lemma}
\newtheorem{prop}[thm]{Proposition}

\theoremstyle{definition}
\newtheorem{defn}{Definition}
\newsymbol\ltimes 226E

\begin{document}

\title{A Groenewold\,-Van Hove Theorem for $S^2$}

\author{Mark J. Gotay}

\address{Department of Mathematics \\University of Hawai`i
\\2565 The Mall \\ Honolulu, HI 96822 USA}

\email{gotay@@math.hawaii.edu}

\author{Hendrik Grundling}

\address{Department of Pure Mathematics \\University of
New South Wales \\P. O. Box 1 \\ Kensington, NSW 2033
Australia}

\email{hendrik@@solution.maths.unsw.edu.au}

\author{C. A. Hurst}

\address{Department of Physics and Mathematical Physics
\\University of Adelaide
\\ G. P. O. Box 498 \\ Adelaide, SA 5001 Australia}

\email{ahurst@@physics.adelaide.edu.au}

\date{February 18, 1995}

\maketitle

\begin{abstract}

We prove that there does not exist a nontrivial
quantization of the Poisson algebra of the symplectic
manifold $S^2$ which is irreducible on the subalgebra
generated by the components $\{S_1,S_2,S_3\}$ of the spin
vector.  We also show that there does not exist such a
quantization of the Poisson subalgebra $\p$ consisting of
polynomials in $\{S_1,S_2,S_3\}$. Furthermore, we show
that the maximal Poisson subalgebra of $\cal P$ containing
$\{1,S_1,S_2,S_3\}$ that can be so quantized is just that
generated by
$\{1,S_1,S_2,S_3\}$.

\end{abstract}

%%%%%%%%%%%%%%%%%%%%%%%%%%%%%%%%%%%%%%%%%%%%%%%%%%%%%%%%%%%%%%%%%%%

\section{Introduction}

In a striking paper, Groenewold \cite{gr} showed that one
cannot ``consistently'' quantize all polynomials in the
classical positions $q^i$ and momenta $p_i$ on $\Bbb
R^{2n}.$ Subsequently Van Hove \cite{vh1,vh2} refined and
extended \gr's result, in effect showing that there does
not exist a quantization functor which is consistent with
the Schr\"odinger quantization of
$\Bbb R^{2n}$. (For discussions of \gr's and \vh's work
and related results, see
\cite{a-m,c,f,go1,g-s,j} and references contained
therein.) However, these theorems rely heavily on certain
properties of $\Bbb R^{2n}$, and so it is not clear
whether they can be generalized. Naturally, one {\em
expects} similar ``no\,-go'' theorems to hold in a wide
range of situations, but we are not aware of any previous
results along these lines.

In this paper we prove a Groenewold-Van Hove theorem for
the symplectic manifold
$S^2$. Our proof is similar to Groenewold's for $\Bbb
R^{2n}$, although it differs from his in several important
respects and is technically more complicated. On the other
hand, as $S^2$ is compact there are no problems with the
completeness of the flows generated by the classical
observables, and so Van Hove's modification of \gr's
theorem is unnecessary in this instance.

To set the stage, let $(M,\omega)$ be a symplectic
manifold. We are interested in quantizing the Poisson
algebra $C^{\infty}(M)$ of smooth real-valued functions on
$M$, or at least some subalgebra $\cal C$ of it, in the
following sense.

\begin{defn} A {\em quantization} of $\cal C$ is a linear
map $\cal Q$ from
$\cal C$ to an algebra of self-adjoint
operators\footnote{Technical difficulties with unbounded
operators will be ignored, as they are not important for
what follows.} on a Hilbert space such that
\ms
\begin{enumerate}\begin{enumerate}
\item[({\em i\/})] $\cal Q\big(\{f,g\}\big) =
-\mbox{i}\big[\cal Q(f),\cal Q(g)\big],$
\end{enumerate}\end{enumerate}

\noindent where $\{\,,\,\}$ denotes the Poisson bracket
and $[\;,\,]$ the commutator. If
$\cal C$ contains the constant function 1, then we also
demand
\ms
\begin{enumerate}\begin{enumerate}
\item[({\em ii\/})] $\cal Q(1) = I.$
\end{enumerate}\end{enumerate}
\end{defn}

As is well known, it is necessary to supplement these
conditions for the quantization to be physically
meaningful. To this end, one often requires that a certain
subalgebra $\cal B$ of observables be represented
irreducibly. Exactly which observables should be taken as
``basic'' in this regard depends upon the particular
example at hand; one typically uses the components of a
momentum map associated to a (transitive) Lie symmetry
group. For $\Bbb R^{2n}$ the relevant group is the
Heisenberg group \cite{f,g-s} and $\cal B = {\text
{span}}\{1,q^i,p_i\,|\, i=1,\ldots,n\}$. In the case of
$S^2$ the appropriate group is SU$(2) \times \Bbb R$
whence the basic observables are span$\{1,S_1,S_2,S_3\}$,
the
$S_i$ being the components of the spin angular momentum.

Alternatively, one could require the {\em strong von
Neumann rule} \cite{vn}

\[\cal Q\big(k(f)\big) = k\big(\cal Q(f)\big)\]

\noindent to hold for all polynomials
$k$ and all $f \in \cal C$ such that $k(f) \in \cal C$.
Usually it is necessary to weaken this condition \cite{f},
insisting only that it hold for a certain subclass $\cal
B$ of observables $f$ and certain polynomials $k$. We
refer to this simply as a ``von Neumann rule.'' In the
case of $\Bbb R^{2n}$, the von Neumann rule as applied to
the $q^i$ and $p_i$ with $k(x) = x^2$ is actually implied
by the irreducibility of the $\cal Q(q^i)$ and $\cal
Q(p_i)$ \cite{c}. But the corresponding statement is not
quite true for $S^2$, as we will see. We refer the reader
to \cite{f} for further discussion of \vn\ rules.

In our view, imposing an irreducibility condition on the
quantization map $\cal Q$ seems more compelling physically
and pleasing aesthetically than requiring
$\cal Q$ to satisfy a \vn\ rule. With this as well as the
observations above in mind, we make

\begin{defn} An {\em admissible quantization} of the pair
$(\cal C,\cal B)$ is a quantization of
$\cal C$ which is irreducible on $\cal B$, where $\cal C
\subset C^\infty(M)$ is a given Poisson algebra, and $\cal
B \subset\cal C$ a given subalgebra.
\end{defn}

The results of \gr\ and Van Hove may then be interpreted
as showing that there does not exist an admissible
quantization of the pair $\big(C^{\infty}(\Bbb
R^{2n}),\,\text{span}\{1,q^i,p_i\,| \linebreak
\, i=1,\ldots,n\}
\big)$ nor, for that matter, of the subalgebra of all
polynomials in the $q^i$ and
$p_i$. We will prove here that there likewise does not
exist an admissible quantization of the pair
$\big(C^{\infty}(S^2),\,\text{span}\{1,S_1,S_2,S_3\}\big)$
nor, for that matter, of the subalgebra of all polynomials
in the components $S_i$ of the spin vector.

\ms

%%%%%%%%%%%%%%%%%%%%%%%%%%%%%%%%%%%%%%%%%%%%%%%%%%%%%%%%%%%%%%%%%%%

\section{No\,-Go Theorems}

Consider a sphere $S^2$ of radius $s > 0$. We view this
sphere as the ``internal'' phase space of a massive
particle with spin $s$ and realize it as the subset of
$\Bbb R^3$ given by \begin{equation} S_1\,\!^2 + S_2\,\!^2
+ S_3\,\!^2 = s^2, \label{cs2} \end{equation}

\noindent where $\bold S = (S_1,S_2,S_3)$ is the spin
angular momentum. The symplectic form is\footnote{Note
that $\omega$ is $1/s$ times the area form $d{\boldsymbol
\sigma}$ on $S^2$. It is the symplectic form on $S^2$
viewed as a coadjoint orbit of SU(2).}

\[\omega =
\frac{1}{2s^2}\sum_{i,j,k=1}^3\epsilon_{ijk}S_i\,dS_j
\wedge dS_k\]

\noindent with corresponding Poisson bracket
\begin{equation}
\{f,g\} = \sum_{i,j,k=1}^3\epsilon_{ijk}S_i\frac{\partial
f}{\partial S_j}
\frac{\partial g}{\partial S_k}
\label{pb}
\end{equation}

\noindent for $f,g \in C^{\infty}(S^2).$ We have the
relations $\{S_i,S_j\} =
\sum_{k=1}^3 \epsilon_{ijk}S_k.$

The group \,SU(2)\, acts transitively on $S^2$ with
momentum map $\bold S= (S_1,S_2,S_3)$, i.e., the pair
$\big(S^2,{\text {SU(2)}}\big)$ is an ``elementary
system'' in the sense of \cite{w}. Thus it is natural to
require that quantization provide an irreducible
representation of \,SU(2). In terms of observables,
quantization should produce a representation which is
irreducible when restricted to the subalgebra generated by
$\{\cx,\cy,\cz\}$. However, this subalgebra does not
include the constants. To remedy this, we consider instead
the central extension ${\text {SU(2)}} \times \Bbb R$ of
\,SU(2)\, by $\Bbb R$ with momentum map $(1,S_1,S_2,S_3)$,
and take the subalgebra generated by these observables to
be the basic set $\cal B$ in the sense of the Introduction.

Let $\cal P$ denote the Poisson algebra of polynomials in
the components
$S_1,S_2,S_3$ of the spin vector $\bold S$ modulo the
relation \eqref{cs2}. (This means we are restricting
polynomials as functions on $\Bbb R^3$ to $S^2$.) We shall
refer to an equivalence class $p \in \cal P$ as a
``polynomial'' and take its degree to be the minimum of
the degrees of its polynomial representatives. We denote
by $\cal P^k$ the subspace of polynomials of degree at
most $k$. In particular, $\cal P^1$ is just the Poisson
subalgebra generated by $\{1,S_1,S_2,S_3\}$.

When equipped with the $L^2$ inner product given by
integration over $S^2$, the vector space
$\p^k$ becomes a real Hilbert space which admits the
orthogonal direct sum decomposition $\cal P^k =
\bigoplus_{l=0}^k \cal H_l$, where $\cal H_l$ is the
vector space of spherical harmonics of degree $l$ (i.e.,
the restrictions to
$S^2$ of homogeneous harmonic polynomials of degree $l$ on
$\Bbb R^3$
\cite{a-b-r}). Note that $\h_1$ is the Poisson subalgebra
generated by
$\{S_1,S_2,S_3\}$. The collection of \shs\
$\big\{Y_l^m,\;l = 0,1,\ldots,k,\;m =
-l,-l+1,\ldots,l\big\}$ forms the standard (complex)
orthogonal basis for the complexification $\p^k_{\Bbb C}$:
\begin{equation*}
\int_{S^2}Y_{l_1}^{m_1*}Y_{l_2}^{m_2}\,d\sigma =
s^2\delta_{l_1l_2}\delta_{m_1m_2}.
\end{equation*}

\noi Thus if $p \in \cal P^k_{\Bbb C}$, we have the
harmonic decomposition
\begin{equation} p = p_k + p_{k-1} + \cdots + p_0,
\label{hd}
\end{equation}

\noi where $p_l \in (\h_l)_{\Bbb C}$ is given by
\begin{equation} p_l =
\frac{1}{s^2}\sum_{m=-l}^l\left(\int_{S^2}Y_l^{m*}p\,d\sigma\right)Y_l^m.
\label{hdy}
\end{equation}

It is well known that O(3) acts orthogonally on $\cal P$,
and that this action is irreducible on each $\h_l$ where
it is the standard real (orbital) angular momentum $l$
representation. The corresponding infinitesimal generators
on
$\h_l$ are $L_i = \{S_i,\cdot\}$. If we identify o(3) and
$\h_1$ as Lie algebras, it follows that the ``adjoint''
action of
$\h_1$ on $\h_l$ given by $S_i \mapsto \{S_i,\cdot\}$ is
irreducible as well.

\ms

Now suppose $\cal Q$ is a quantization of $\cal P$, so
that \begin{equation}
\big[\cal Q(S_i),\cal Q(S_j)\big] =
\mbox{i}\sum_{k=1}^3\epsilon_{ijk}\cal Q(S_k) \label{com}
\end{equation}
\noi and
\vspace{-1ex}
\begin{equation}
\cal Q\big(\bold S^2\big) = s^2I.
\label{s2}
\end{equation}
\vspace{1ex}
\noindent If in addition $\q$ is admissible on $(\cal
P,\,\cal P^1)$, then
$\p^1$ must be irreducibly represented. Since as a Lie
algebra $\cal P^1$ is isomorphic to \,su$(2) \times \Bbb
R$, its irreducible representations are all
finite-dimensional.\footnote {Every irreducible
representation of su$(2) \times \Bbb R$  by (essentially)
self-adjoint operators on an invariant dense domain in a
Hilbert space can be integrated to a continuous
irreducible representation of $\text {SU(2)}
\times \Bbb R$ \cite[\S 11.10.7.3]{b-r}. But it is well
known that every such representation of this group is
finite-dimensional.} These are just the usual (spin
angular momentum) representations labeled by
$j = 0,\frac{1}{2},1,\ldots$, where\footnote{In what
follows we use standard quantum mechanical notation, cf.
\cite{m}.}
\begin{equation}
\sum_{i=1}^3\cal Q(S_i)^2 = j(j+1)I.
\label{j2}
\end{equation}

\noindent The Hilbert space corresponding to the quantum
number $j$ has dimension
$2j+1$, with the standard orthonormal basis $\big\{|\,j,
m\rangle,\: m = -j,-j+1,\ldots,j\big\}$ consisting of
eigenvectors of $\cal Q(S_3)$. We regard the
representation defined by $j=0$ as trivial, in that it
corresponds to quantum spin 0.

The following result shows that admissibility implies a
weak type of \vn\ rule on $\cal P^1$.

\begin{prop} If $\cal Q$ is an admissible quantization of
$(\cal P,\,\cal P^1)$, then
\begin{equation}
\cal Q\big(S_i\,\!^2\big) = a \cal Q(S_i)^2 + cI \label{ii}
\end{equation}

\noindent for $i = 1,2,3$, where $a$ and $c$ are real
constants with $a^2 + c^2
\neq 0$. Furthermore, for $i \neq \ell$,
\begin{equation}
\cal Q(S_iS_\ell) = \frac{a}{2}\big(\cal Q(S_i)\cal
Q(S_\ell)+\cal Q(S_\ell)\cal Q(S_i)\big).
\label{pr}
\end{equation}
\label{qvn}
\end{prop}
\vspace{-2ex} We have placed the proof, which is rather
long and technical, in Appendix A so as not to interrupt
the exposition.

Observe that on summing \eqref{ii} over $i$, we get
$s^2=aj(j+1)+3c$ which fixes the constant $c$ in terms of
$s,\:j$ and $a.$

{}From these relations the main result now follows.

\begin{thm}[No\,-Go Theorem] There does not exist a
nontrivial admissible quantization of $(\cal P,\cal P^1)$.
\label{ng2}
\end{thm}
\begin{pf} Suppose there did exist an admissible
quantization of $(\cal P,\cal P^1)$; we shall show that
for $j>0$ this leads to a contradiction.

First observe that we have the classical equality

\[s^2S_3 = \{S_1\,\!^2 - S_2\,\!^2,S_1S_2\} -
\{S_2S_3,S_3S_1\}.\]

\noindent Quantizing this, a calculation using \eqref{ii},
\eqref{pr},
\eqref{com} and \eqref{j2} gives

\[s^2\qz = a^2\bigg (j(j+1) - \frac {3}{4} \bigg ) \cal
Q(S_3).\]

\noindent Thus either $\cal Q(S_3) = 0$, whence $j=0$, or
$j > 0$, in which case
\begin{equation} s^2 = a^2\bigg (j(j+1) - \frac {3}{4}
\bigg ). \label{qc1}
\end{equation}

\noindent Observe that when $j = \frac{1}{2}$, this
implies that $s = 0$, which is impossible. Henceforth take
$j>\frac{1}{2}.$

Next we quantize the relation

\[2s^2\cy\cz = \big\{\cy\,\!^2,\{\cx\cy,\cx\cz\}\big\} -
\frac{3}{4}\big\{\cx\,\!^2,\{\cx\,\!^2,\cy\cz\}\big\}.\]

\noi Using \eqref{pr}, \eqref{ii}, \eqref{com} and
\eqref{j2}, the l.h.s. becomes

\[as^2\big(\q(S_2)\q(S_3)+\q(S_3)\q(S_2)\big) = as^2 \big
(2\qy\qz -
\mbox{i}\qx\big )\]

\noi while the r.h.s. reduces to

\[a^3\bigg (j(j+1) - \frac{9}{4}\bigg) \big(2\qy\qz -
\mbox{i}\qx\big ).\]

\noindent Since for $j > \frac{1}{2}$ the matrix element

\[\big{\langle}\, j, j\,|\,2\qy\qz - \mbox{i}\qx|\,j,
j\!-\!1\big{\rangle} =
\text {i} \bigg (\frac{1}{2} - j\bigg )\sqrt{2j}\]

\noindent is nonzero, it follows that

\[as^2 = a^3\bigg (j(j+1) - \frac{9}{4} \bigg ).\]

\noindent If $a=0$ \eqref{qc1} yields $s=0$, whereas if $a
\neq 0$ this conflicts with \eqref{qc1}. Thus we have
derived contradictions provided $j > 0$. Since $j=0$ is
the trivial representation, the theorem is proven. \end{pf}

This contradiction shows that the quantization goes awry
on the level of quadratic polynomials. On the other hand,
there are many admissible quantizations of the Poisson
subalgebra $\cal P^1$ of all polynomials of degree at most
one, viz. the irreducible representations $\q$ of \,su$(2)
\times \Bbb R$ with $\q(1) = I$. Thus it is of interest to
determine the largest subalgebra of $\cal P$ containing
$\{1,S_1,S_2,S_3\}$ that can be so quantized. We will now
show that this largest subalgebra is just $\cal P^1$
itself. Unfortunately, this is not entirely
straightforward, since $\cal P^1$ is not a maximal Poisson
subalgebra of $\cal P$; indeed, if $\cal O$ denotes the
Poisson subalgebra of odd polynomials (i.e., polynomials
all of whose terms are of odd degree), then
$\cal P^1$ is contained in $\tilde{\cal O} = \cal O \oplus
\Bbb R$.

To prove the result, we proceed in two stages. First we
show that $\cal P^1$ is maximal in $\tilde {\cal O}$, and
then prove a no\,-go theorem for
$\tilde {\cal O}$.

\begin{prop}
$\cal P^1$ is maximal in $\tilde{\cal O}$. \label{max}
\end{prop}

\begin{pf} Actually, the constants are unimportant, and it
will suffice to prove that the Poisson subalgebra $\h_1$
generated by $\{S_1,S_2,S_3\}$ is maximal in
$\cal O$.

Set $\cal O^l = \cal O \cap
\p^l$, where henceforth $l$ is odd. For $k$ odd it is
clear from \eqref{pb} and
\eqref{cs2} that $\{\h_k,\h_l\} \subset \oo^{k+l-1}$.

Let $\cal R$ be the Poisson algebra generated by a single
polynomial $r \in
\cal O^l$ of degree $l>1$ together with $\h_1$. Evidently
$\cal R \subset \cal O$; we must show that $\cal O \subset
\cal R$. We will accomplish this in a series of lemmas.

\begin{lem} If in its harmonic decomposition an element of
$\cal R$ has a nonzero component in $\h_k$, then
$\h_k\subset\cal R$. \label{inhomo}
\end{lem}

\begin{pf} Let $\cal R' \subset \cal R$ be the span of all
elements of the form
\[\big\{h_n,\ldots\big\{h_2,\{h_1,r\}\big\}\ldots\big\}\]

\noi for $h_i \in \h_1$ and $n \in \Bbb N.$ Then $\cal R'$
is an o(3)-invariant subspace of $\oo^l \subset \p$. Since
the representation of o(3) on
$\cal P$ is completely reducible, $\cal R'$ must be the
direct sum of certain
$\h_k$ with $k \leq l.$ Consequently, if when harmonically
decomposed an element of $\cal R'$ has a nonzero component
in some $\h_k$, then $\h_k \subset \cal R'
\subset \cal R$.
\renewcommand{\qedsymbol}{$\bigtriangledown$} \end{pf}

Now by assumption $r_l \neq 0$ in $\h_l$ and hence $\h_l
\subset \cal R.$ Then
$\{\h_l,\h_l\} \subset
\oo^{2l-1} \cap \cal R$. We will use this fact to show
that $\oo^{2l-1}
\subset \cal R$. The proof devolves upon an explicit
computation of the harmonic decomposition of
$\{Y_l^m,Y_l^n\}$.

\begin{lem} For each $j$ in the range $0 < j \leq 2l$, we
have \begin{equation*}
\{Y_l^{l-j},Y_l^l\} =
\sum_{k=1}^{l}y_{2k-1}(l-j,l)Y_{2k-1}^{2l-j}.
\end{equation*} In particular, when $j = 1$ the top
coefficients $y_{2l-1}(l-1,l)$ are nonzero. Furthermore,
provided $l\geq 5$, $k \geq \frac{l-1}{2}$ and $k > l -
\frac{j+1}{2}$, the coefficients $y_{2k-1}(l-j,l)$ are
nonzero.
\label{main}
\end{lem} Since the proof requires an extended
calculation, we defer it until Appendix B.

\begin{lem}
$\oo^{2l-1} \subset \cal R$.
\end{lem}

\begin{pf} Decompose $Y_l^m = R_l^m+\text {i}I_l^m$ into
real and imaginary parts, with
$R_l^m,\;I_l^m\in\cal H_l$. So in $\cal P_{\Bbb C}$ we
observe that both
$\Re\{Y_l^m,\,Y_l^n\}=\{R_l^m,R_l^n\}- \{I_l^m,I_l^n\}$ and
$\Im\{Y_l^m,Y_l^n\}=\{R_l^m,I_l^n\}+\{I_l^m,R_l^n\}$
belong to
$\{\cal{H}_l,\,\cal{H}_l\}$.

Thus if the harmonic
decomposition of
$\{Y_l^m,Y_l^n\}$ has a nonzero $k^{\text{th}}$ component,
then either its real or imaginary part must be nonzero
which allows us to conclude that
$\{\cal{H}_l,\cal{H}_l\}$ contains an element with nonzero
component in
$\cal{H}_k$.

If $l=3$, then $\h_3 \subset \cal R$. Now consider the
bracket
$\{Y_3^{2},\,Y_3^3\}$. By Lemma \ref{main} with $j=1$ it
has a nonzero
$5^{\text{th}}$ component, so by the preceding and Lemma
\ref{inhomo} it follows that $\h_5 \subset \cal R$. Since
by definition
$\h_1\subset\cal R$, we then have $\oo^5 = \h_1 \oplus
\h_3 \oplus \h_5 \subset
\cal R.$

If $l\geq 5$, we consider
$\{Y_l^{-2},\,Y_l^l\}\in\cal R_{\Bbb C}$. By Lemma
\ref{main} with $j=l+2$, the preceding and Lemma
\ref{inhomo} we conclude that

\[\h_{l-2}\oplus\h_l\oplus\cdots\oplus\h_{2l-1}\subset\cal
R.\]

\noi Hence $\h_{l-2}\subset\cal R$, so by the same
argument applied to
$\{Y_{l-2}^{-2},Y_{l-2}^{l-2}\}$ we get that
$\h_{l-4}\subset\cal R$. Continuing in this way we obtain
$\h_{l-2n}\subset\cal R$ for all $n$ with $l-2n\geq 3$. In
particular, taking
$n = \frac{l-3}{2}$, we get
$\h_3\subset\cal R$. But we have already remarked that
$\h_1 \subset \cal R,$ so the lemma is proven.
\renewcommand{\qedsymbol}{$\bigtriangledown$} \end{pf}

Thus $\cal R$ must contain all odd polynomials of degree
at most $2l-1$. To obtain higher degree polynomials, we
need only bracket $\h_{2l-1} \subset \cal R$ with itself
and apply the argument above to conclude that $\cal
O^{4l-3}
\subset \cal R$. Continuing in this manner, we have
finally that $\cal O \subset
\cal R$, and this proves the proposition. \end{pf}

Our strategy in proving the no\,-go theorem for $\tilde
{\oo}$ is the same as for $\p$. To begin, we use
admissibility to obtain a weak version of a cubic
\vn\ rule.

\begin{prop} If $\q$ is an admissible quantization of
$(\tilde{\cal O},\p^1)$, then \begin{equation}
\q\big(S_i\,\!^3\big) = a \qi ^3 + c \qi
\label{vn3}
\end{equation}

\noi for $i=1,2,3,$ where $a$ and $c$ are real constants.
Furthermore, when $i
\neq \ell$,
\begin{equation}
\q(\ci\cj\ci) = a\qi\qj\qi + \frac{1}{3}(a+c)\qj.
\label{iji}
\end{equation}
\noi Finally,
\begin{equation}
\q(\cx\cy\cz) = a\qx\qy\qz + \frac{a}{2 \mathrm{i}}\big
(\qx^2 - \qy^2 +
\qz^2\big).
\label{123}
\end{equation}
\label{cvn}
\end{prop}
\vspace{-1.5ex} Again the proof is placed in Appendix A.

We derive some consequences of these results. Multiplying
\eqref{cs2} through by $S_\ell$ and quantizing gives
\[\sum_{i=1}^3 \cal Q(S_iS_\ell S_i) = s^2 \cal
Q(S_\ell).\]

\noindent Applying \eqref{iji} and \eqref{vn3}, this in
turn becomes
\begin{equation} a\sum_{i=1}^3 \cal Q(S_i)\cal
Q(S_\ell)\cal Q(S_i) = \bigg(s^2 - \frac{2a}{3} -
\frac{5c}{3}\bigg) \cal Q(S_\ell). \label{s3s}
\end{equation}

\noi We also find, by rearranging the factors in
$\sum_{i=1}^3 \cal Q(S_i)\cal Q(S_\ell)\cal Q(S_i)$, that
\begin{equation}
\sum_{i=1}^3 \cal Q(S_i)\cal Q(S_\ell)\cal Q(S_i) = \big
(j(j+1) -1\big )\qj.
\label{s3j}
\end{equation}

\noi A comparison of \eqref{s3s} and \eqref{s3j} yields
\begin{equation} s^2 = a\bigg(j(j+1) -\frac{1}{3}\bigg ) +
\frac{5c}{3} \label{js}
\end{equation}

\noi provided $j>0.$

\begin{thm} There does not exist a nontrivial admissible
quantization of
$(\tilde{\cal O},\p^1)$. \label{ng3}
\end{thm}
\vspace{-3ex}
\begin{pf} Suppose $\q$ were an admissible quantization of
$(\tilde{\cal O},\p^1)$; we will show that then $j=0.$

Consider the Poisson bracket relation
\begin{eqnarray*} 3s^4S_3 &\!\!\! =\!\!\! &
4\{S_1\,^3,S_2S_3\,^2\} - 4\{S_2\,^3,S_3\,^2S_1\} +
\{S_2\,^2S_1,S_2\,^3\} \\  & &
\rule{0ex}{3ex} -
\,\{S_2S_1\,^2,S_1\,^3\} - 6\{S_2\,^3,S_1\,^3\}
-3\{S_2S_3\,^2,S_3\,^2S_1\}.
\end{eqnarray*}

\noindent Upon quantizing, an enormous calculation using
\eqref{vn3},
\eqref{iji}, \eqref{com} and \eqref{s3j}
gives\footnote{This calculation was done using the {\sl
Mathematica} package {\sl NCAlgebra} \cite{h-m}.}
\begin{eqnarray*} 3s^4 \cal Q(S_3) & = & \bigg(3a^2j^4 +
6a^2j^3 + 14acj^2 + 8a^2j^2 \nonumber \\ & & \mbox{} +
14acj + 5a^2j + \frac{29c^2}{3} - 14
\frac{ac}{3} - \frac{7a^2}{3}\bigg)\cal Q(S_3) \nonumber
\\ \rule{0ex}{3ex} & &
\mbox{} - (10a^2 + 4ac)\cal Q(S_3)^3 \end{eqnarray*}

\noindent which in view of \eqref{js} simplifies to
\begin{equation}
\hspace{-.25in}
\bigg [\frac{4}{3}(2a-c)(a+c) - a(7a+4c)j(j+1)\bigg ]\cal
Q(S_3) + a(10a+4c)\qz^3 = 0.
\label{cc1}
\end{equation}

Now suppose $j = \frac{1}{2}$, so that $\qz^2 =
\frac{1}{4}I$. Then
\eqref{cc1} implies that $a = -4c$ which, when substituted
into \eqref{js} yields $s=0$. Similarly, when $j = 1$,
$\qz^3 = \qz$. In this case \eqref{cc1} implies that $a =
-c$, and again \eqref{js} requires $s=0.$  Thus we have
derived contradictions for these two values of $j$.
Henceforth take $j>1.$

Next we quantize

\[ 6s^2\cx\cy\cz = \{S_1\,\!^3,\cyy\cx\}
+\{S_2\,\!^3,\czz\cy\} +
\{\czzz,\cxx\cz\}. \]

\noi Another computer calculation using \eqref{vn3},
\eqref{iji}, \eqref{123},
\eqref{com} and \eqref{s3j} yields
\begin{eqnarray*}
\lefteqn{6s^2\Big[a\qx\qy\qz +
\frac{a}{2\text{i}}\big(\qx^2 - \qy^2 + \qz^2\big)\Big]}
\nonumber \\
\rule{0ex}{1ex} & &
\hspace{-.18in} \rule{0ex}{4ex} =
-3a\text{i}\big(c-2a+aj(j+1)\big)\big[\qx^2-\qy^2+\qz^2 +
2\text{i}\qx\qy\qz\big].
\label{cc2a}
\end{eqnarray*}

\noi Since for $j>1$ the matrix element
\begin{eqnarray*}
\lefteqn{\hspace{-1.5in} \Big {\langle}\,j,
j\!-\!2\,\Big|\,\qx\qy\qz +
\frac{1}{2\text{i}}\big(\qx^2 -
\qy^2 +
\qz^2\big)\Big|\,j, j\Big{\rangle}} \\ & & =
\rule{0ex}{4ex}\frac{1}{2{\text {i}}}(1-j)\sqrt{j(2j-1)}
\nonumber
\end{eqnarray*}
\noi is nonzero, we conclude that either $a=0$ or
\begin{equation} s^2 = c-2a + aj(j+1).
\label{cc2}
\end{equation}
\noi If $a=0$, \eqref{cc1} implies that $c=0$, and then
\eqref{js} leads to a contradiction, so
\eqref{cc2} must hold. Subtracting \eqref{cc2} from
\eqref{js} gives $c = -5a/2$; substituting this into
\eqref{cc1} produces

\[3a^2(j^2+j-3)\qz = 0.\]

\noi But then $j = (-1 \pm \surd{13})/2$, neither of which
is permissible.

Thus we have derived contradictions for all $j>0$ and the
theorem follows.
\end{pf}

We remark that Theorem \ref{ng3} is actually sharper than
Theorem \ref{ng2}; we have included the latter because it
is simpler.

As Proposition \ref{max} shows, by augmenting $\cal P^1$
with a single odd polynomial, we generate $\tilde{\cal
O}$. Similarly we can generate all of $\cal P$ from $\cal
P^1$ and a single polynomial not in
$\tilde{\cal O}$, which implies that the only Poisson
subalgebras of $\cal P$ strictly containing
$\cal P^1$ are $\tilde{\cal O}$ and $\cal P$ itself. This
can be proven in the same manner as Proposition \ref{max},
but in the interests of economy, we make do with:

\begin{lem} Any Poisson subalgebra of $\cal P$ which
strictly contains $\cal P^1$ also contains $\tilde{\cal
O}$. \end{lem}

\begin{pf} By Proposition \ref{max} it suffices to
consider the Poisson algebra
$\cal T$ generated by $\cal P^1$ and a polynomial $p$ not
in
$\tilde{\cal O}$. Then $p$ has a component in some
$\h_{2k}$ for $k>0$, so by Lemma \ref{inhomo} it follows
that $\h_{2k}\subset\cal T$. Now consider the bracket
$\{Y_{2k}^{2k-1},\, Y_{2k}^{2k}\}$. According to Lemma
\ref{main}, either its real or imaginary part has a
nonzero component in
$\h_{4k-1}$. Hence $\h_{4k-1}\subset\cal T$, and so by
Proposition \ref{max} we have $\tilde{\cal O}\subset\cal
T$.
\end{pf}

Now given any Poisson subalgebra of $\cal P$ strictly
containing $\cal P^1$ we only have to apply Theorem
\ref{ng3} to the subalgebra $\tilde{\cal O}$ inside it to
obtain a contradiction, hence:

\begin{thm} No nontrivial quantization of $\cal P^1$ can
be admissibly extended beyond $\cal P^1$.
\label{dogo}
\end{thm}

This result stands in marked contrast to the analogous one
for $\Bbb R^{2n}.$ There one runs into difficulties with
{\em cubic} polynomials in
$\{1,q^i,p_i\},$ so that $\cal P^2$ (the Poisson algebra
of polynomials of degree at most two) is a maximal
polynomial subalgebra containing $\cal P^1$ that can be
admissibly quantized \cite{g-s}. This dichotomy seems to
be connected with the fact that for $\Bbb R^{2n}$, $\cal
P^2$ is the Poisson normalizer of $\cal P^1,$ whereas for
$S^2$ the normalizer of $\cal P^1$ is itself. On the other
hand, it should be noted that there are other maximal
polynomial subalgebras of
$C^{\infty}(\Bbb R^{2n})$ containing $\cal P^1$ that can
be admissibly quantized, for instance the Schr\"odinger
subalgebra

\[\left\{\sum_{i=1}^nh^i(q^1,\dots,q^n)p_i +
k(q^1,\dots,q^n)\right\}\]

\noi where the $h^i$ and $k$ are polynomials. Here one
encounters problems when  one tries to extend to terms
which are quadratic in the momenta.

\ms

Finally, a word is in order regarding the case $j = 0$ --
the one instance in which we did not derive a
contradiction. It happens that the spin 0 representation
of $\p^1$ {\em can\/} be extended, in a unique way, to an
admissible quantization
$\q$ of all of $\p$. Indeed, given $p \in \p$, let $p_0$
denote the constant term in the harmonic decomposition
\eqref{hd} of $p$. Then $\q:\p \rightarrow
\Bbb C$ defined by $\q(p) = p_0$ is, technically, an
admissible quantization of
$(\p,\p^1)$. To prove this, it is only necessary to show
that $\q$ so defined is a Lie algebra homomorphism, which
in this context means $\{p,p'\}_0 = 0$. But from
\eqref{hdy} and \eqref{pb}, in vector notation,

\[4\pi s^2\{p,p'\}_0 = \int_{S^2}\{p,p'\}\,d\sigma =
\int_{S^2}\bold S
\cdot (\nabla p \times \nabla p')\,d\sigma = s\int_{S^2}
(\nabla p \times \nabla p') \cdot d{\boldsymbol \sigma},
\]

\noi which vanishes by the divergence
theorem.\footnote{This actually is a consequence of a
general fact about momentum maps on compact symplectic
manifolds, cf. \cite[p.~ 187]{g-s}.}

To show that this quantization is unique, we need

\begin{lem} For $l > 0$, $\h_l = \{\h_1,\h_l\}.$
\label{l1l}
\end{lem}

\begin{pf} For $l>0$ $\{\h_1,\h_l\}$ is a nontrivial
invariant subspace of
$\h_l$, and hence by the irreducibility of the
$\h_1$-action on $\h_l$,
$\{\h_1,\h_l\} = \h_l.$
\end{pf}

Now suppose $\chi$ is an admissible quantization of
$(\p,\p^1)$ with $j=0$ so that by the above, $\chi$ is a
linear map $\p \ra \Bbb C$ which must annihilate Poisson
brackets. But then $\chi(p) = \chi(p_0)$, since by Lemma
\ref{l1l} each term $p_l \in
\h_l$ for $l > 0$ in the harmonic decomposition of $p$
must be a sum of terms of the form $\{h_l,r_l\}$ for some
$h_l \in \h_1$ and $r_l \in
\h_l$. It follows that $\chi$ is uniquely determined by
its value on the constants, and as
$\chi(1) = I = \q(1)$ we must have $\chi = \q.$

Although the corresponding representation of $\p^1$ is
trivial in that $\qi = 0$ for all $i$, it is worth
emphasizing that $\q$ is {\em not\/} zero on the remainder
of $\p$. For example,
$\q(\cii) =
\frac{s^2}{3}I$ for all $i$, consistent with Proposition
\ref{qvn} and
\eqref{s2}. The existence of this trivial yet ``not
completely trivial'' representation of $\p$ -- for any
nonzero value of the classical spin $s$ -- may be related
to a well-known ``anomaly'' in the geometric quantization
of spin, cf. \cite[\S7]{t} and \cite[\S11.2]{s}.

\ms\ms\ms

%%%%%%%%%%%%%%%%%%%%%%%%%%%%%%%%%%%%%%%%%%%%%%%%%%%%%%%%%%%%%%%%%%%

\section{Discussion}

Theorems \ref{ng2} and \ref{dogo} could have been
``predicted'' on the basis of geometric quantization
theory. Here one knows that one can quantize those
classical observables $f \in C^{\infty}(M)$ whose flows
preserve a given polarization $P$. In our case we take $P$
to be the antiholomorphic polarization on $S^2$ (thought
of as $\Bbb CP^1$); then $\p^1$ is exactly the set of
polarization-preserving observables. However, in geometric
quantization theory one does not expect to be able to
consistently quantize observables outside this class
\cite{w}.

Further corroboration for our results is provided by
Rieffel \cite{r}, who showed that there are no strict
SU(2)-invariant deformation quantizations of
$C^{\infty}(S^2)$. In fact, it seems that only the
polynomial algebra $\cal P^1
\subset C^{\infty}(S^2)$ can be rigorously deformation
quantized in an SU(2)-invariant way \cite{k}.

There are several points we would like to make concerning
the \vn\ rules for
$S^2$, especially in comparison with those for $\Bbb
R^{2n}$. Our no\,-go theorems may be interpreted as
stating that the ``Poisson bracket $\ra$ commutator'' rule
is {\em totally} incompatible with even the relatively weak
\vn\ rules given in Propositions \ref{qvn} and \ref{cvn}.
On $\Bbb R^{2n}$, on the other hand, the \vn\ rule $k(x) =
x^n$ with $n=2$ does hold for $x\in\p^1
\subset \p^2$ in the metaplectic representation
\cite{c,g-s}, and with any
$n \geq 0$ for $x=q^i$ in the Schr\"odinger
representation. Moreover it is curious that, according to
Propositions
\ref{qvn} and \ref{cvn}, requiring that the $S_i$ be
irreducibly represented does not yield strict \vn\ rules
as happens for $\{q^i,p_i\}$ on $\Bbb R^{2n},$ cf.
\cite{gr,c}.\footnote{In the case of $\Bbb R^{2n}$, one
usually employs the Stone-\vn\ theorem to show that
irreducibility implies the \vn\ rules, but in fact one can
prove this without recourse to the Stone-\vn\ theorem.} It
is not clear why $S^2$ and $\Bbb R^{2n}$ behave
differently in these regards.

We remark that it is substantially easier to prove the
no\,-go theorems \ref{ng2} and \ref{ng3} if one assumed
from the start that the strict \vn\ rules $\q(\cii) =
\q(\ci)^2$ and $\q(\ciii) = \q(\ci)^3$ hold. We can gain
some insight into this as follows. Suppose the strong \vn\
rule applied to one of the $S_i$, say
$S_3$, so that $\cal Q(S_3\,\!^n) = \cal Q(S_3)^n$ for all
positive integers
$n$. Suppose furthermore that $\q$ is injective when
restricted to the polynomial algebra $\Bbb R[S_3]$
generated by $S_3$. Provided it is appropriately
continuous,
$\q$ will then extend to an isomorphism from the real
$C^*$--algebra $C^*(S_3)$ (consisting of the closure of
$\Bbb R[S_3]$ in the supremum norm on $S^2$ with pointwise
operations) to the real $C^*$--algebra generated by $\cal
Q(S_3)$. But this implies that the classical spectrum of
$S_3$ (i.e., the set of all values it takes) is the same
as the operator spectrum of $\cal Q(S_3)$. Since the
classical spectrum of $S_3$ is $[-s,s]$ whereas the
quantum spectrum is discrete, it is clear -- in retrospect
-- why no (strong) \vn\ rule can apply to
$S_3$.

Since $S^2$ is in a sense the opposite extreme from $\Bbb
R^{2n}$ insofar as symplectic manifolds go, our result
lends support to the contention that no\,-go theorems
should hold in some generality. Nonetheless, these two
examples are special in that they are symplectic
homogeneous spaces \big($\Bbb R^{2n}$ for the translation
group $\Bbb R^{2n}$, and $S^2$ for $\mbox{SU}(2)$\big).
Thus in both cases we are quantizing finite-dimensional
Lie algebras \big(the Heisenberg algebra and $\mbox{su}(2)
\times \Bbb R$, respectively, which are certain central
extensions by $\Bbb R$ of $\Bbb R^{2n}$ and su(2)\big).
Will a similar analysis work for other symplectic
homogeneous spaces, e.g., $\Bbb CP^n$ with group
$\mbox{SU}(n+1)$? How does one proceed in the case of
symplectic manifolds which do not have such a high degree
of symmetry? What set of observables will play the role of
the distinguished subalgebras generated by
$\{1,q^i,p_i\}$ and $\{1,S_1,S_2,S_3\}$ (which are the
components of the momentum mappings for the Hamiltonian
actions of the Heisenberg group and
$\mbox{SU}(2) \times \Bbb R$, resp.)? For a cotangent
bundle $T^*Q$, the obvious counterpart would be the
infinite-dimensional abgebra of linear momentum
observables $P_X + f,$ where $P_X$ is the momentum in the
direction of the vector field $X$ on $Q$ and $f$ is a
function on $Q$. (These are the components of the momentum
mapping for the transitive action of
$\text{Diff}(Q) \ltimes C^{\infty}(Q)$ on $T^*Q$.) In this
regard, it is known that geometric quantization (formally)
obeys the von Neumann rules $\cal Q\big(P_X\,^2\big) =
\cal Q(P_X)^2$ and $\cal Q(f^n) = \cal Q(f)^n$ \cite{go2}.
We hope to explore some of these issues in future papers.

\ms\ms

%%%%%%%%%%%%%%%%%%%%%%%%%%%%%%%%%%%%%%%%%%%%%%%%%%%%%%%%%%%%%%%%%%%

\section*{Acknowledgments}

One of us (M.J.G.) would like to thank L. Bos, G. Emch, G.
Goldin, M. Karasev and J. Tolar for enlightening
conversations, and the University of New South Wales for
support while this research was underway. He would also
like to express his appreciation to both the organizers of
the XIII Workshop on Geometric Methods in Physics and the
University of Warsaw, who provided a lively and congenial
atmosphere for working on this problem. This research was
supported in part by NSF grant DMS-9222241.

Another one of us (H.G.) would like to acknowledge the
support of an ARC--grant.

\ms

%%%%%%%%%%%%%%%%%%%%%%%%%%%%%%%%%%%%%%%%%%%%%%%%%%%%%%%%%%%%%%%%%%%

\appendix

\section{\vn\ Rules}

Here we provide the proofs of Propositions \ref{qvn} and
\ref{cvn}.

Fix an irreducible representation of $\text {su(2)} \times
\Bbb R$ labeled by the quantum number $j$. The proofs
depend on the following three facts:
\vskip 2pt
\begin{enumerate}
\item[({\em i\/})] As the representation is irreducible,
any endomorphism of the representation space is an element
of the enveloping algebra of the generators $\qi$
\cite[Prop. 2.6.5]{d}, and hence can be expressed as a
polynomial in the $\qi$.

\vskip 6pt

\item[({\em ii\/})] Under the su(2)-action, a monomial
$\qx\,\!^{n_1}\qy\,\!^{n_2}\qz\,\!^{n_3}$ of degree $|n| =
n_1 + n_2 + n_3$ transforms as a tensor operator of rank
$|n|$.

\vskip 6pt

\item[({\em iii\/})] Under the induced action of su(2) on
$\p$, a monomial $S_1\,\!^{n_1}S_2\,\!^{n_2}S_3\,\!^{n_3}$
transforms as a symmetric tensor of rank $|n|$. Since the
quantization $\q$ is (infinitesimally) equivariant, it
follows that
$\q(S_1\,\!^{n_1}S_2\,\!^{n_2}S_3\,\!^{n_3})$ also
transforms as a symmetric tensor operator of rank $|n|$.
\end{enumerate}

\noi When $ |n| > 1$, the tensor operators
$\qx\,\!^{n_1}\qy\,\!^{n_2}\qz\,\!^{n_3}$ and
$\q(S_1\,\!^{n_1}S_2\,\!^{n_2}S_3\,\!^{n_3})$ are
reducible.

Equation \eqref{ii} then follows from the observation that
as $\qii$ is a reducible symmetric tensor operator of rank
2, its irreducible constituents must be of even rank, and
hence it must be equal to an even polynomial in the $\qi$
of degree at most 2. Equations
\eqref{pr}, \eqref{vn3}, \eqref{iji} and \eqref{123}
follow from similar observations. The rest of the argument
is the working out of these observations.

Letting $n = (n_1,n_2,n_3)$ be a multi-index of length
$|n| = n_1+n_2+n_3$, we denote
$\qx\,\!^{n_1}\qy\,\!^{n_2}\qz\,\!^{n_3}$ by $\Bbb Q_{\bf
n}^{|n|}$ and $S_1\,\!^{n_1}S_2\,\!^{n_2}S_3\,\!^{n_3}$ by
$\Bbb S_{n}^{|n|}$. Since the commutation relations
\eqref{com} are nonlinear, the difference between each
$\Bbb Q_{\bf n}^{|n|}$ and its symmetrization $\Bbb Q_{\bf
(n)}^{|n|}$ is a linear combination of tensor operators
$\Bbb Q_{\bf m}^{|m|}$ of lower rank $|m|$. Thus we may
use the symmetrized tensor operators $\Bbb Q_{\bf
(n)}^{|n|}$ as a basis for the enveloping algebra of the
operators $\qi$. Then by ({\em i\/}) we can expand
\begin{equation}
\q\big(\sm\big) = \sum_{|n|=0}^{d}[\,m\,|\,n\,]\,\Bbb
Q_{\bf (n)}^{|n|}
\label{Expand}
\end{equation}

\noi as a polynomial of degree $d$, say.\footnote{It can
be shown that $d
\leq 4j$, but we do not need this fact.} We can further
decompose each
 \begin{equation}
\Bbb Q_{\bf (n)}^{|n|} = \sum_{\ld = 0}^{|n|}\sum_{\mu =
-\ld}^{\ld}(\,n\,|\,\ld\;\mu\,)\,T^{\ld}_{\mu}
\label{t}
\end{equation}

\noi into a sum of irreducible spherical tensor operators
$T^{\ld}_{\mu}$ of rank $\ld$ (with $\mu =
-\ld,\dots,\ld$),\footnote {Since they are symmetric the
tensor operators $\Bbb Q_{\bf (n)}^{|n|}$  are ``simply
reducible,'' i.e., in the decomposition \eqref{t} there is
at most one irreducible constituent $T^{\ld}_{\mu}$ for
each weight $\ld$. This follows from a consideration of
Young tableaux. Likewise, the
$\q\big(\sm\big)$, being symmetric, are simply reducible,
and hence there is no degeneracy in either \eqref{Irredc}
or \eqref{Thirdord} below.} so that
\eqref{Expand} becomes
\begin{equation}
\q\big(\sm\big) =
\sum_{|n|=0}^{d}\sum_{\lambda=0}^{|n|}\sum_
{\mu=-\lambda}^{\lambda}[\,m\,|\,n\,]\,(\,n\,|\,\lambda\;
\mu\,)\,T_{\mu}^{\lambda}. \label{Exdecom} \end{equation}

On the other hand, for $|m|=2$ we may directly decompose
\begin{equation}
\q\big(\Bbb S^2_m\big) = \sum_{\nu=-2}^2
(\,m\,|\,2\;\nu\,)\,V_{\nu}^2 + (\,m\,|\,0\;0\,)\,V^0_0,
\label{Irredc}
\end{equation}

\noi where the irreducible constituents $V^{\eta}_{\nu}$
of $\q\big(\Bbb S^2_m\big)$ are given by
\begin{eqnarray} V_{\nu}^2 & =
&\sum_{|m|=2}(\,2\;\nu\,|\,m\,)\,\q(\Bbb S^2_m),
\label{Irreda}\\ V_0^0 & =
&\sum_{|m|=2}(\,0\;0\,|\,m\,)\,\q(\Bbb S^2_m).
\label{Irredb}
\end{eqnarray}

\noi Here we have used the relation
$\sum_{\ld \geq 0}\sum_{\nu =
-\ld}^{\ld}(\,m\,|\,\ld\;\nu\,)\,(\,\ld\;\nu\,|\,m'\,) =
\delta_{mm'}$. Note that, as it is symmetric, $\q\big(\Bbb
S^2_m\big)$ has no irreducible rank 1 constituent.
Combining \eqref{Irreda} and \eqref{Irredb} with
\eqref{Exdecom}, we obtain
\begin{eqnarray} V_{\nu}^2 & =
&\sum_{|m|=2}\sum_{|n|=0}^{d}\sum_{\lambda =0}^
{|n|}\sum_{\mu =
-\ld}^{\ld}(\,2\;\nu\,|\,m\,)\,[\,m\,|\,n\,]\,(\,n\,|\,\lambda
\;\mu\,)\,T_{\mu}^{\lambda}, \label{Exirreda}\\  V^0_0 & =
&\sum_{|m|=2}\sum_{|n|=0}^{d}\sum_{\lambda=0}^{|n|}
\sum_{\mu=-\ld}^{\ld}(\,0\;0\,|\,m\,)\,[\,m\,|\,n\,]\,(\,n\,|\,\lambda\;\mu\,)\,
T_{\mu}^{\lambda}.
\label{Exirredb}
\end{eqnarray}

Now apply a rotation $R \in {\text {SU(2)}}$ to
\eqref{Exirreda} to obtain
\begin{eqnarray*}
\lefteqn{\sum_{\nu'=-2}^2D_{\nu\nu'}^2(R)V_{\nu'}^2 =} \\
& &
\sum_{|m|=2}\sum_{|n|=0}^{d}\sum_
{\lambda=0}^{|n|}\sum_{\mu=-\ld}^{\ld}\sum_{\mu'=-\ld}^{\ld}(\,2\;\nu\,|\,m\,)\,
[\,m\,|\,n\,]\, (\,n\,|\,\lambda\;\mu\,)\,D_{\mu
\mu'}^{\lambda}(R)\,T_{\mu'}^{\lambda}.
\label{Rotate}
\end{eqnarray*}

\noi Multiplying both sides of this by
$D_{\nu\rho'}^2(R)^*$ and integrating over the group
manifold, the orthogonality theorem for products of
representations \cite[\S 7.1]{b-r} yields \begin{equation}
V_{\rho'}^2 = \left[\sum_{|m|=2}\sum_{|n|=0}^{d}
(\,2\;\nu\,|\,m\,)\,[\,m\,|\,n\,]\,(\,n\,|\,2\;\nu\,)\right]T_{\rho'}^2.
\label{FConda}
\end{equation}

\noi Applying the same procedure to \eqref{Exirredb} but
using
$D^0_{00}(R)^*$ instead, we similarly obtain
\begin{equation}  V_0^0 =
\left[\sum_{|m|=2}\sum_{|n|=0}^{d}(\,0\;0\,|\,m\,)\,[\,m\,|\,n\,]\,
(\,n\,|\,0\;0\,)\right]T_0^0.\label{FCondd}
\end{equation}

{}From \eqref{FConda} and \eqref{FCondd}, we see that
$V_{\rho'}^2$ is proportional to $T_{\rho'}^2$, and
$V_0^0$ to $T^0_0$; let the corresponding constants of
proportionality be $a$ and $b$. Substituting
\eqref{FConda} and \eqref{FCondd} into \eqref{Irredc} and
inverting
\eqref{t} gives
\begin{eqnarray}
\q\big(\Bbb S^2_m\big) & = &
a\sum_{\nu=-2}^2(\,m\,|\,2\;\nu\,)\,T_{\nu}^2 +
b(\,m\,|\,0\;0\,)\,T_0^{0} \nonumber \\  & = &
a\sum_{|n|=2}\sum_{\nu=-2}^2(\,m\,|\,2\;\nu\,)\,(\,2\;\nu\,|\,n\,)\,\Bbb
Q^2_{(n)} + b(\,m\,|\,0\;0\,)\,(\,0\;0\,|\,0\,)\,\Bbb Q^0_0
\nonumber \\ \rule{0ex}{4ex} & = &
a\sum_{|n|=2}^{\,}P_{2}(m,n)\,\Bbb Q^2_{(n)} + bP_0(m,0)I.
\label{Secord} \end{eqnarray}
\noi In this last expression,
$P_{\ld}(m,n) =
\sum_{\nu=-\ld}^{\ld}(\,m\,|\,\ld\;\nu)\,(\,\ld\;\nu\,|\,n\,)$
is the projector which picks off the $m^{\text
{th}}$-component (with respect to
$\Bbb Q^{|m|}_{(m)}$) of the irreducible rank
$\ld$ constituent of
$\Bbb Q_{(n)}^{|n|}$. Since $m$ has length 2, $m$ must be
of the form $1_i+1_\ell$, where $1_i$ is the multi-index
with 1 in the
$i^{\text{th}}$-slot and zeros elsewhere. Similarly, when
$|n| =2$, $n = 1_p+1_q$. Then we have
\begin{eqnarray*}  P_{2}(1_i+1_\ell,1_p+1_q) & = &
\frac{1}{2}\big(\delta_{ip}\delta_{\ell
q}+\delta_{iq}\delta_{\ell p}\big) -
\frac{1}{3}\delta_{i\ell}\delta_{pq},
\label{Proj2} \\ \rule{0ex}{4ex} P_0(1_i+1_\ell,0) & = &
\frac{1}{3}\delta_{i\ell}.
\label{Proj0}
\end{eqnarray*}

Setting $\ell = i$, \eqref{Secord} and \eqref{j2} then
yield
\begin{equation*}
\qii=a\left(\qi^2-\frac{1}{3}j(j+1)I\right)+\frac{1}{3}bI.
\label{Prop2}
\end{equation*}

\noi The constants $a$ and $b$ must both be real \big(as
$\qii$ is self-adjoint\big), and both cannot be
simultaneously zero \big(for this would contradict
\eqref{s2}\big). Thus, upon setting $c =
\big(b-aj(j+1)\big)/3$, we obtain \eqref{ii}.

Taking $\ell \neq i$, we similarly obtain \eqref{pr}.

The same arguments can be applied, {\em mutatis mutandis},
to $\q\big(\Bbb S^3_m\big)$. Then \eqref{Secord} is
replaced by
\begin{eqnarray}
\q\big(\Bbb S^3_m\big) & = &
a\sum_{\nu=-3}^3(\,m\,|\,3\;\nu\,)\,T_{\nu}^3 +
b\sum_{\nu=-1}^1(\,m\,|\,1\;\nu\,)T_{\nu}^1 \nonumber\\  &
= & \rule{0ex}{4ex} a\sum_{|n|=3}P_{3}(m,n)\,\Bbb
Q^3_{(n)} + b\sum_{|n|=1}P_{1}(m,n)\,\Bbb Q^1_n.
\label{Thirdord}
\end{eqnarray}

\noi The projectors in this instance are
\begin{eqnarray*}
\lefteqn{\hspace
{-1.17in}P_{3}(1_i+1_k+1_{\ell},1_p+1_q+1_r)} \\ & = &
\rule{0ex}{4ex}{1\over
6}\sum\delta_{ip}\delta_{kq}\delta_{\ell r}- {1\over
30}\sum\delta_{ik}(\delta_{pq}\delta_{\ell
r}+\delta_{pr}\delta_{\ell q}+\delta_{qr}\delta_{\ell p}),
\label{eq:Proj3}\\ \rule{0ex}{4ex}
P_{1}(1_i+1_k+1_{\ell},1_p) & = & {1\over
15}(\delta_{ik}\delta_{\ell
p}+\delta_{k\ell}\delta_{ip}+\delta_{\ell i}\delta_{kp}),
\label{eq:Proj1} \end{eqnarray*}

\noi where the sums are over all permutations of
$\{i,k,\ell\}$.

Setting $i = k = \ell$, \eqref{Thirdord} reduces to
\begin{equation*}
\qiii=a\left(\qi^3-\frac{3}{5}\Big(j(j+1)-\frac{1}{3}\Big)\qi\right)+\frac{1
}{5}b\qi.
\label{Prop3}
\end{equation*}

\noi Again by self-adjointness, $a$ and $b$ must be real.
Upon setting $c =
\frac{1}{5}\big(b-3\big[j(j+1) -
\frac{1}{3}\big]\big)$ and rearranging the above, we obtain
\eqref{vn3}.

Finally, taking $i=k \neq \ell$ and $(i,k,\ell) = (1,2,3)$
in
\eqref{Thirdord}, we similarly obtain
\eqref{iji} and \eqref{123}, respectively.

Although it is not apparent from this derivation, it can be
shown that \eqref{pr} is actually a consequence of
\eqref{ii} and, likewise, that \eqref{iji} and \eqref{123}
both follow from \eqref{vn3}.

\ms\ms

%%%%%%%%%%%%%%%%%%%%%%%%%%%%%%%%%%%%%%%%%%%%%%%%%%%%%%%%%%%%%%%%%%%

\appendix
\addtocounter{section}{1}
\section{On the Harmonic Decomposition of
$\{Y_l^{l-j},Y_l^l\}$}

In this appendix we prove Lemma \ref{main}, which computes
the harmonic decomposition of $\{Y_l^{l-j},Y_l^l\}$ for
$l>0$. Specifically, for each $j,\; 0 < j \leq 2l$, we have
\begin{equation*}
\{Y_l^{l-j},Y_l^l\} =
\sum_{k=1}^{l}y_{2k-1}(l-j,l)Y_{2k-1}^{2l-j},
\end{equation*}

\noi where the coefficients $y_{2k-1}(l-j,l)$ have the
following properties:

\begin{enumerate}
\item[({\em i\/})] when $j = 1,$ the top coefficients
$y_{2l-1}(l-1,l) \neq 0$, and
\ms
\item[({\em ii\/})] provided $l\geq 5$, $k \geq
\frac{l-1}{2}$ and $k > l -
\frac{j+1}{2}$, we have
$y_{2k-1}(l-j,l) \neq 0.$
\end{enumerate}
\ms The proof will be presented in several steps. We refer
the reader to
\cite[Chapter XIII and Appendix C]{m} for the relevant
background and conventions on spherical harmonics.

%%%%%%%%%%%%%%%%%%%%%%%%

\ms

{\em Step 1.} It is convenient to work in spherical
coordinates $(\theta,\phi)$ on $S^2$. The generators $L_3$
and $L_{\pm} = L_1 \pm \text {i}L_2$ of O(3) then take the
form

\[L_3 = \frac{\partial}{\partial \phi} \;\;\mbox{and}\;\;
L_{\pm} = \pm\mbox{i} e^{\pm \text{i}\phi} \bigg
(\frac{\partial}{\partial \theta} \pm \mbox{i}\cot
\theta \frac{\partial}{\partial \phi}\bigg). \]

\noi They satisfy

\[L_3Y_l^m = {\text i}mY_l^m, \;\;\;\;L_+Y_l^m = {\text
i}\beta_{l,m}Y_l^{m+1}\;\;\mbox{ and }\;\; L_-Y_l^m =
{\text i}\beta_{l,m-1}Y_l^{m-1}\]

\noi where $\beta_{l,m} = \sqrt{(l+m+1)(l-m)}.$

The Poisson bracket \eqref{pb} becomes
\begin{equation*}
\{f,g\} = \frac{\csc \theta}{s} \bigg (\frac{\partial
f}{\partial
\phi}\frac{\partial g}{\partial \theta} - \frac{\partial
f}{\partial
\theta}\frac{\partial g}{\partial \phi}\bigg ),
\end{equation*}

\noi which can be rewritten in terms of the $L_i$ as
\begin{equation*}
\{f,g\} = -\frac{\text{i}}{s}\csc \theta\,
e^{-{\text{i}}\phi}\Big (L_3(f)L_+(g)-L_+(f)L_3(g)\Big).
\label{pbs} \end{equation*}

\ms

%%%%%%%%%%%%%%%%%%%%%%%%%%%%%%%%%%%%%%%%

{\em Step 2.} As a specific instance of this formula, we
compute \begin{equation}
\{Y_l^m,Y_l^n\}\sin \theta \, e^{\text{i}\phi} =
\frac{\mbox
i}{s}\Big(m\,\beta_{l,n}Y_l^mY_l^{n+1}-n\,\beta_{l,m}Y_l^{m+1}Y_l^n\Big).
\label{mn}
\end{equation}

\noi Now we have the product decomposition \begin{eqnarray}
Y_{l_1}^{m_1}Y_{l_2}^{m_2}&=&\sum_{l = |l_1-l_2|}^{l_1+l_2}
\sqrt{\frac{(2l_1+1)(2l_2+1)}{4\pi(2l+1)}}\:
\langle\,l_1\,l_2\,0\,0\,|\,l\,0\,\rangle \nonumber\\ &
&\times\:\langle\,l_1\,l_2\,m_1\,m_2\,|\,l\;m_1\!+\!m_2\,\rangle
Y_l^{m_1+m_2},
\label{pd}
\end{eqnarray}

\noi where the quantities
$\langle\,l_1\,l_2\,m_1\,m_2\,|\,L\,M\,\rangle$ are \cg s.
Applying \eqref{pd} to the r.h.s. of \eqref{mn} gives
\begin{eqnarray}
\{Y_l^m,Y_l^n\}\sin \theta \, e^{\text{i}\phi} & = &
\frac{\mbox i}{s}
\sum_{j=0}^{2l}\frac{2l+1}{\sqrt{4\pi(2j+1)}}\langle\,l\,l\,0\,0\,|\,j\,0\,
\rangle
\label{yyr} \\ & & \mbox{}\!\! \times
\Big
[m\,\beta_{l,n}\langle\,l\,l\,m,\,n\!+\!1\,|\,j,\,m\!+\!n\!+1\,\rangle
\nonumber \\ & & \mbox{} \;\;\;\;\; - n\, \beta_{l,m}
\langle\,l\,l\,m\!+\!1,\,n\,|\,j,\,m\!+\!n\!+1\,\rangle
\Big ]Y_j^{m+n+1}.
\nonumber
\end{eqnarray}

On the other hand, since $\{\h_l,\h_l\} \subset
\oo^{2l-1}$ we can expand
\begin{equation*} \{Y_l^m,Y_l^n\} =
\sum_{k=1}^{l}\sum_{r=-2k+1}^{2k-1}y_{2k-1,r}(m,n)Y_{2k-1}^r.
\end{equation*}

\noi Multiply this equation through by $\sin
\theta\,e^{\text {i}\phi}$. Since
$\sin \theta\,e^{\text {i}\phi} = -
\sqrt{8\pi/3}\,Y_1^1$, we can apply the product
decomposition \eqref{pd} to the r.h.s. of the resulting
expression to obtain \begin{eqnarray}
\{Y_l^m,Y_l^n\}\sin \theta \, e^{\text{i}\phi} & = & -
\sqrt{\frac{8\pi}{3}}
\bigg
(\sum_{k=1}^{l}\sum_{r=-2k+1}^{2k-1}y_{2k-1,r}(m,n)\big
[Y_{2k-1}^rY_1^1\big ]\bigg)
\label{yyl} \\ & = &
\sum_{k=1}^{l}\sum_{r=-2k+1}^{2k-1}y_{2k-1,r}(m,n)\Big
(e_{2k-1,r}Y_{2k-2}^{r+1} - d_{2k-1,r}Y_{2k}^{r+1}\Big),
\nonumber \end{eqnarray}

\noi where
\[d_{t,r} =
\sqrt{\frac{(t+r+1)(t+r+2)}{(2t+1)(2t+3)}}\:\:\: \mbox{
and }\:\:e_{t,r} =
\sqrt{\frac{(t-r)(t-r-1)}{(2t+1)(2t-1)}}.\]

Comparing \eqref{yyr} with \eqref{yyl}, we see that $r =
m+n$ and hence
\begin{eqnarray*}
\lefteqn{\hspace{-.5in} \sum_{k=1}^{l}y_{2k-1,m+n}(m,n)
\Big (e_{2k-1,m+n}Y_{2k-2}^{m+n+1} -
d_{2k-1,m+n}Y_{2k}^{m+n+1}\Big)} \\ & = &
\frac{\mbox i}{s}
\sum_{j=0}^{2l}\frac{2l+1}{\sqrt{4\pi(2j+1)}}\langle\,l\,l\,0\,0\,|\,j\,0\,
\rangle
\\ & &
\mbox{} \!\times
\Big
[m\,\beta_{l,n}\langle\,l\,l\,m,\,n\!+\!1\,|\,j,\,m\!+\!n\!+1\,\rangle
\\ & &\mbox{}\!\;\;\;\;\;\; - n\, \beta_{l,m}
\langle\,l\,l\,m\!+\!1,\,n\,|\,j,\,m\!+\!n\!+1\,\rangle
\Big ]Y_j^{m+n+1}.
\end{eqnarray*}

\noi Note that the sum on the l.h.s. contains only even
degree harmonics, and hence the sum on the r.h.s. must as
well. \big(This is reflected by the vanishing of the \cg s
$\langle\,l\,l\,0\,0\,|\,j\,0\,\rangle$ for $j$ odd, cf.
\cite[eq. (C23.a)]{m}.\big) Thus we may reindex $j=2k,\: k
= 0,\ldots,l$ on the r.h.s.. Upon reindexing $k \mapsto
k+1$ in the first term on the l.h.s. and equating
coefficients of $Y_{2k}^{m+n+1}$ on both sides, we obtain
the recursion relation
\begin{eqnarray} y_{2k-1}(m,n) & = &
\frac{e_{2k+1,m+n}}{d_{2k-1,m+n}}\,y_{2k+1}(m,n) \nonumber
\\ & & \mbox{} -
\frac{\text{i}}{s} \frac{2l+1}{\sqrt{4\pi(4k+1)}\,
d_{2k-1,m+n}}\, \langle
\,l\,l\,0\,0\,|\,2k\,0\, \rangle \nonumber \\ & & \mbox{}
\times \Big{[}
m\,\beta_{l,n}\langle\,l\,l\,m\,n\!+\!1\,|\,2k,\,m\!+\!n\!+1\,\rangle
\nonumber
\\ & & \mbox{}\;\;\;\;\;\; - n\, \beta_{l,m}
\langle\,l\,l\,m\!+\!1,\,n\,|\,2k,\,m\!+\!n\!+1\,\rangle
\Big{]} \label{rry}
\end{eqnarray}

\noi where we have abbreviated $y_{2k-1,m+n}(m,n) =:
y_{2k-1}(m,n)$.

\ms

{\em Step 3.} Now we specialize even further, setting $m =
l-j$ and $n=l$ for $j = 1,\cdots,2l$. Then we have
\begin{equation}
\{Y_l^{l-j},Y_l^l\} =
\sum_{k=1}^{l}y_{2k-1}(l-j,l)Y_{2k-1}^{2l-j}, \label{ylyl}
\end{equation}

\noi as in the statement of the lemma, and \eqref{rry}
reduces to
\begin{eqnarray} y_{2k-1}(l-j,l) & = &
\frac{e_{2k+1,2l-j}}{d_{2k-1,2l-j}}\,y_{2k+1}(l-j,l)
\nonumber \\ & & \mbox{} +
\frac{\text{i}}{s} \frac{l(2l+1)}{\sqrt{4\pi(4k+1)}\,
d_{2k-1,2l-j}}\,\beta_{l,l-j}
\langle \,l\,l\,0\,0\,|\,2k\,0\, \rangle \nonumber \\ & &
\mbox{} \times
\langle\,l\,l\,l\!-\!j\!+\!1,\,l\,|\,2k,\,2l\!-\!j\!+1\,\rangle.
\label{rrl}
\end{eqnarray}

\noi The main reason for this choice of $m$ and $n$ is
that it simplifies the last term in \eqref{rry}, since
$\beta_{l,l}=0$.

Before proceeding, we must evaluate the \cg s in
\eqref{rrl}. Using the Racah formula \cite[(C.21), (C.22)
and (C.23b)]{m}, we compute

\[\langle \,l\,l\,0\,0\,|\,2k\,0 \,\rangle =
(-1)^{k+l}\sqrt{4k+1}\,
\sqrt{\frac{(2l-2k)!}{(2l+2k+1)!}}\,
\frac{(2k)!(l+k)!}{(k!)^2(l-k)!}\]

\noi and
\begin{eqnarray*}
\lefteqn{\hspace{-1in}
\langle\,l\,l\,l\!-\!j\!+\!1,\,l\,|\,2k,\,2l\!-\!j\!+1\,\rangle
= } \hspace{-1in}
\\ & & \\ & & \sqrt{4k+1}\,
\sqrt{\frac{(2l)!(j-1)!(2k+2l-j+1)!}{(2l-j+1)!(2l+2k+1)!(2l-2k)!(2k-2l+j-1)!}}.
\end{eqnarray*}

\noi Substituting these expressions as well as those for
$d_{2k-1,2l-j}$,
$e_{2k+1,2l-j}$ and $\beta_{l,l-j}$ into \eqref{rrl}, the
recursion relation becomes
\begin{eqnarray*}
\tilde y_{2k-1} & = &
\sqrt{\frac{(2k-2l+j+1)(2k-2l+j)(4k-1)}
{(2k+2l-j)(2k+2l-j+1)(4k+3)}}\,\tilde y_{2k+1} \\ & &
\mbox{} + (-1)^k
\sqrt{\frac{(2k+2l-j-1)!}{(2k-2l+j-1)!}}\,
\frac{\sqrt{4k-1}(4k+1)(l+k)!(2k)!}{(2l+2k+1)!(l-k)!(k!)^2},
\end{eqnarray*}

\noi where $\tilde y_{2k-1}$ is defined according to

\[y_{2k-1}(l-j,l) = \frac{\text i}{s}(-1)^l
\frac{l(2l+1)}{\sqrt{4\pi}}\,\sqrt{\frac{j!(2l)!}{(2l-j)!}}\,\tilde
y_{2k-1}.\]

\noi Since $l>0$, $y_{2k-1}(l-j,l) = 0$ iff $\tilde
y_{2k-1} = 0$.

Finally, we rewrite this in the form
\begin{equation}
\tilde y_{2k-1} = Z_{2k-1}\big [\tilde y_{2k+1} + (-1)^k
W_{2k-1}\big ]
\label{rrb}
\end{equation}

\noi where

\[Z_{2k-1}= \sqrt{\frac{(2k-2l+j+1)(2k-2l+j)(4k-1)}
{(2k+2l-j)(2k+2l-j+1)(4k+3)}}\]

\noi and

\[W_{2k-1} = \sqrt{\frac{(2k+2l-j+1)!}{(2k-2l+j+1)!}}\,
\frac{\sqrt{4k+3}(4k+1)(l+k)!(2k)!}{(2l+2k+1)!(l-k)!(k!)^2}.
\]

\noi Notice that for fixed $j$ and $l$, $y_{2k-1}(l-j,l) =
0$ if $2k-1 \leq 2l-j-2.$

Upon setting $k=l$ and $j=1$, \eqref{rrb} gives $\tilde
y_{2l-1} = (-1)^l Z_{2l-1}W_{2l-1} \neq 0$ and hence
$y_{2l-1}(l-1,l) \neq 0$. This proves ({\em i\/}).

%%%%%%%%%%%%%%%%%%%%%%%

\ms

{\em Step 4.} Fix $l$ and $j$, in which case the
coefficients $Z_{2k-1}$ and
$W_{2k-1}$ are nonzero whenever $2k-1 > 2l-j-2$. We claim
that \begin{equation} Z_{2k-1}W_{2k-1} < W_{2k-3}
\label{ratio}
\end{equation}

\noi provided $l \geq 5$ and $k \geq \frac{l-1}{2}$.
Indeed, we compute

\[\frac{Z_{2k-1}W_{2k-1}}{W_{2k-3}} =
\frac{(4k+1)(l-k+1)(2k-1)}{(4k-3)(2l+2k+1)k}.\]

\noi But now one verifies that the maximum of the r.h.s.
of this expression is 18/25 on the domain in the
$k,l$-plane determined by the above inequalities along
with the fact that $k \leq l.$

%%%%%%%%%%%%%%%%%%%%%%%

\ms

{\em Step 5.} We can now prove statement ({\em ii\/}) of
the lemma. Fix $l\geq 5$ and
$j\in\{1,\ldots,2l\}$. The solution of the recursion
relation \eqref{rrb} with initial condition $\tilde
y_{2l+1} = 0$ is

\[\tilde y_{2l-2n-1} = (-)^l \sum_{k=0}^n(-1)^{k}U_k\]

\noi where $U_k=
W_{2l-2k-1}\prod\limits_{t=0}^{n-k}Z_{2l-2n+2t-1}$, and
$n=0,1,\ldots,\text{max}\{\frac{l+1}{2},\frac{j-1}{2}\}$.
For the alternating sum $\sum_{k=0}^n(-1)^kU_k$ we have by
\eqref{ratio} that $U_{k+1}/U_k>1$, so if
$n$ is even the sum is

\[U_0+(U_2-U_1)+\cdots+(U_n-U_{n-1})>0\]

\noi since $U_0$ and all bracketed terms are positive, and
if $n$ is odd

\[(U_0-U_1)+(U_2-U_3)+\cdots+(U_{n-1}-U_n)<0\]

\noi likewise. So the sum is always nonzero in the given
ranges of the parameters, hence $\tilde y_{2k-1}\not=0$ as
claimed.

%%%%%%%%%%%%%%%%%%%%%%%%%%%%%%%%%%%%%%%%%%%%%%%%%%%%%%%%%%%%%%%%%%%

\ms\ms

%%%%%%%%%%%%%%%%%%%%%%%%%%%%%%%%%%%%%%%%%%%%%%%%%%%%%%%%%%%%%%%%%%%

\end{document}